\documentstyle[11pt]{article}

\newcommand{\vra}{Vitale, R.A. }

\newcommand{\norm}{\|}
\newtheorem{theorem}{\bf Theorem}

\newtheorem{rmk}{\bf Remark}

\newtheorem{ex}{\bf Exercise}

\def\thinsquare{\vcenter{\hrule height 0.3pt \hbox{\vrule width0.3pt height6pt
\kern6.5pt \vrule width 0.3pt}
\hrule height 0.3pt}}
\def\qed{\hfill$\thinsquare$}

\newcommand{\be}{\begin{equation}}
\newcommand{\ee}{\end{equation}}

\newcommand{\proof}{\noindent {\bf Proof  }}

\newcommand{\clconv}{\overline{\rm conv}}

\newcommand {\ZZ}{{Z\kern-.45em Z}}
\newcommand {\RR}{{I\kern-.3em R}}

\newcommand {\CC}{{I\kern-.6em C}}
\newcommand {\NN}{{I\kern-.3em N}}

\renewcommand{\phi}{\varphi}

\newcommand{\reciptwopi}{\frac{1}{\sqrt{2\pi}}}
\newcommand{\dist}{{\delta}}
\newcommand{\zd}{Z^{(d)}}
\newcommand{\mor}{M/\dist}
\newcommand{\td}{\tilde{d}}
\begin{document}


\begin{center}
{\large \bf INTRINSIC $L_p$ METRICS FOR CONVEX BODIES}
\vspace{.3cm}

{\bf RICHARD A. VITALE}

\end{center}

\vspace{.5cm}

\begin{quote}
{\bf Abstract.} Intrinsic $L_p$ metrics are defined and shown to satisfy a dimension--free 
bound with respect to the Hausdorff metric. \newline MSC 2000: 52A20, 52A27, 52A40, 60G15.
\end{quote}

\pagestyle{empty}

\section{Introduction} $L_p$ metrics for convex bodies were originally introduced in the context of approximation questions (\cite{v75_mcl}; also \cite{arn89}, \cite{bor1}, \cite{bor2}, \cite{lud}) and later shown to be comparable to the Hausdorff metric:

\be \label{first}
\delta_p(K,L)\geq	 c_{p,d,K,L}\cdot \delta^{(p+d-1)/p}(K,L)
\ee
(\cite[Corollary 1]{v85_Lp}; also, \cite[Lemma 1]{grosch}).  Bounds like (\ref{first}) have been useful for establishing stability results and others ([2], [3], [6]--[21], [25]--[27], [29]).

The dimension $d$ of the underlying space appears not only in the form of (\ref{first}) but also in the definition of $\delta_p$ itself. A natural question to ask is whether these dependencies can be avoided. Analogous to the renormalization of quermassintegrals to intrinsic volumes (\cite{mcm75, sch93}), this amounts to asking whether there are dimension--free, or {\em intrinsic}, versions of the $L_p$ metrics and of (\ref{first}). A positive answer was given in \cite{v93_class}. In this note, we give an improved version of that result and some related comments. 

Let $K,L$ be convex bodies in $\RR^d$, and let $Z^{(d)}=(Z_1,Z_2, \ldots, Z_d)$ be a vector of independent, standard Gaussian variables. For $1\leq p < \infty$, the {\it intrinsic} $L_p$ {\it metric} is given by

\be \label{intdefin}
\dist_p^{*}(K,L)=\left[c_pE|h_K(\zd)-h_L(\zd)|^p\right]^{1/p},
\ee
where $c_p=1/E|Z_1|^p={\pi}^{1/2}\left[2^{p/2}\Gamma\left((p+1)/2\right)\right]^{-1}$ is chosen so that ${\dist_p^*}(\{x\}, \{\tilde{x}\})=\norm x - \tilde{x}\norm$ for any $x, \tilde{x} \in \RR^d$ (cf. \cite[Eqn. 22]{v93_class}). This coincides with the usual $L_p$ metric up to a multiplicative constant (which depends on both $p$ and $d$). We begin with an explicit proof of the following:

\begin{theorem} ${\dist_p^*}$ is intrinsic, $1\leq p < \infty$.
\end{theorem}

\proof Suppose that $K,L$ lie in a proper subspace of $\RR^{d}$: without loss of generality, $K,L \subset \mbox{span}\{(x_1,x_2,\ldots,
x_{\td}, \underbrace{0,0,\ldots,0}_{d-\td})\}$. Let $\sigma : \RR^d \rightarrow \RR^{\td}$ be the associated projection operator. For any $x\in \RR^d$, one has $h_K(x)=h_{\sigma K}(x)=h_K(\sigma x),$
and the explicit form of $E|h_K(\zd)-h_L(\zd)|^p$ gives

\begin{eqnarray*}
&\phantom{=}&\int\limits_{z_1=-\infty}^{\infty}\cdots\int\limits_{z_d=-\infty}^{\infty}|h_K(z)-h_L(z)|^p 
(2\pi)^{-d/2}\Pi_{i=1}^de^{-z_i^2/2}dz_d \cdots dz_1\\
&=&\int\limits_{z_1=-\infty}^{\infty}\cdots\int\limits_{z_d=-\infty}^{\infty}|h_K(\sigma z)-h_L(\sigma z)|^p 
(2\pi)^{-d/2}\Pi_{i=1}^de^{-z_i^2/2}dz_d\cdots dz_1\\
&=&\int\limits_{z_1=-\infty}^{\infty}\cdots\int\limits_{z_{\td}=-\infty}^{\infty}|h_K(\sigma z)-h_L(\sigma z)|^p 
(2\pi)^{-{\td}/2}\Pi_{i=1}^{\td}e^{-z_i^2/2}dz_{\td}\cdots dz_1\\
&=&E|h_{\sigma K}(Z^{(\td)})-h_{\sigma L}(Z^{(\td)})|^p,
\end{eqnarray*}
so that ${\dist_p^*}(K,L)={\dist_p^*}(\sigma K, \sigma L).$ \qed 

\section{An Intrinsic Bound} 
We now give an intrinsic form of (\ref{first}).
 \begin{theorem}
For $p\geq 1$ and finite dimensional convex bodies $K,L$:
\be \label{bound}
{\dist_p^*}(K,L)\geq (1/4) \dist(K,L) e^{-\frac{1}{2\pi}\left(\frac{V_1(\clconv (K\cup L))}{\dist(K,L)}\right)^2}.
\ee
\end{theorem}

\proof For notational convenience, let $\dist=\dist (K,L)$ and $V_1=V_1(\clconv (K\cup L)).$ Referring to \cite[Theorem 1]{v93_class}, 
let $M$ be the unique, positive solution to
\be \label{mdef}
E\left(M-\dist Z\right)_+ =\reciptwopi V_1\;,
\ee
which also satisfies $E|h_K(\zd )- h_L(\zd )|^p \geq E\left[\left(\dist Z - M\right)_+\right]^p.$ From (\ref{intdefin}), it follows that $${\dist_p^*}(K,L) \geq \left( c_p E\left[\left(\dist Z - M\right)_+\right]^p \right)^{1/p}.$$ Now
\begin{eqnarray*}
E\left[\left(\dist Z - M\right)_+\right]^p&=& \dist^p E\left[\left( Z-\mor \right)_+ \right]^p \\
&=& \dist^p \int_{\mor}^{\infty} \left(z-\mor \right)^p\reciptwopi e^{-z^2/2}dz \\
&=&\dist^p \int_0^{\infty} y^p \reciptwopi e^{-(y+\mor )^2/2}dy\\
&\geq&\dist^p \int_0^{\infty} y^p \reciptwopi e^{-y^2-(\mor )^2}dy \\
&\geq&\dist^p e^{-(\mor )^2}\int_0^{\infty} y^p \reciptwopi  e^{-y^2}dy \\
&\geq&\dist^p e^{-(\mor )^2}2^{-(p+1)/2}\int_0^{\infty} w^p \reciptwopi e^{-w^2/2}dw \\
&\geq&\dist^p e^{-(\mor )^2}2^{-(p+1)/2}\frac{1}{2c_p}. 
\end{eqnarray*}
Therefore,
\be
{\dist_p^*}(K,L) \geq \left[c_p \dist^p e^{-(\mor )^2}2^{-(p+1)/2}\frac{1}{2c_p}\right]^{1/p}\geq (1/4) \dist e^{-(\mor )^2}. \label{penult}
\ee
From (\ref{mdef}), one has $M=E\left(M - \dist Z \right)\leq E\left(M- \dist Z \right)_+ = \reciptwopi V_1$, which can then be substituted into (\ref{penult}). \qed

\bigskip 
\noindent {\bf Remarks}
\smallskip


\noindent 1. The interested reader may want to compare (\ref{bound}) with \cite[eqn. 24]{v93_class}.
\smallskip

\noindent 2. Theorem \ref{bound} has an equivalent formulation for Gaussian processes. Suppose that $K,L$ are convex bodies in Hilbert space and that $\{X_t\}_{t\in K}$, $\{X_t\}_{t \in L}$ are corresponding isonormally indexed, mean--zero, bounded Gaussian processes. Then 

$$E|\sup_{t\in K} X_t - \sup_{t\in L}X_t|\geq (1/4) \delta
e^{-B^2/{\delta^2}},$$
where $\delta$ is the Hausdorff distance between $K$ and $L$, and $B=$ \linebreak $E\max\{\sup_{t\in K}X_t, \sup_{t\in L}X_t\}.$

\smallskip

\noindent 3. It is possible to extend Theorem 2 to all so-called {\it GB} convex bodies in Hilbert space. In this case, $Z^{(d)}$ is replaced in (\ref{bound}) by $Z^{(\infty)}=(Z_1,Z_2,\dots)$, an infinite sequence of independent standard Gaussian variables. One can ask then for the metric space completion of the class of convex bodies in Hilbert space under $\dist_p^*$. Unfortunately this turns out to have limited geometric significance. This can be seen using some facts from Gaussian processes: let $\{e_n\}_n$ be an orthonormal basis and $a_n=(\log (n+1))^{-1/2}$. Define $K_N=\clconv \{a_ne_n\}_1^{N}.$ For any $p$, this is a Cauchy sequence, and the limit (in the completion) can be identified with $\clconv\{a_ne_n\}_1^{\infty}$. Now let $\tilde{K}_N=\clconv\{a_ne_n\}_N^{\infty}.$ Each of these is also in the completion. Moreover, for any $p$, they form a Cauchy sequence  whose limit is almost surely a (strictly) positive constant. But this cannot be a supremum $h_K(Z^{(\infty)})$ of Gaussian random variables for any $K$. Thus the completion goes beyond the natural geometric setting. An alternate approach is given in \cite{v00_intvols}.


\hfill Department of Statistics, U-4120

\hfill University of Connecticut 

\hfill Storrs, CT 06269-4120 USA

\vspace{.2cm}

\hfill r.vitale@uconn.edu
\end{document}